\documentclass[10pt]{article}
\pdfoutput=1

\usepackage{authblk}
\usepackage{bm}
\usepackage[hidelinks]{hyperref}
\usepackage{url}
\usepackage{amssymb,amsthm}
\usepackage{amsmath}
\usepackage[noabbrev]{cleveref}
\usepackage{cite}
\usepackage{glossaries}

\newcommand*{\glsargt}{t}
\renewcommand*{\glsentryfmt}{%
  \let\orgglsarg\glsargt
  \ifdefempty\glsinsert
  {}%
  {%
    \let\glsargt\glsinsert
    \let\glsinsert\relax
  }%
  \glsgenentryfmt
  \let\glsargt\orgglsarg
}
\newcommand*{\glsargi}{i}
\renewcommand*{\glsentryfmt}{%
  \let\orgglsarg\glsargi
  \ifdefempty\glsinsert
  {}%
  {%
    \let\glsargi\glsinsert
    \let\glsinsert\relax
  }%
  \glsgenentryfmt
  \let\glsargi\orgglsarg
}

\makeglossaries
\newglossaryentry{t}{name=\ensuremath{t},description={time},text={\ensuremath{{\glsargt}}}}
\newglossaryentry{T}{name=\ensuremath{T},description={optimization horizon}}

\newglossaryentry{x}{name=\ensuremath{\boldsymbol{x}_t},description={state at time \gls{t}}, user1={\ensuremath{\boldsymbol{x}}}}
\newglossaryentry{u}{name=\ensuremath{\boldsymbol{u}_t},description={input at time \gls{t}},user1={\ensuremath{\boldsymbol{u}}}}
\newglossaryentry{d}{name=\ensuremath{\boldsymbol{d}_t},description={disturbance at time \gls{t}},user1={\ensuremath{\boldsymbol{d}}}}
\newglossaryentry{yh}{name=\ensuremath{\boldsymbol{z}_{t-h:t}},description={output history at time \gls{t}},user1={\ensuremath{\boldsymbol{z}}}}
\newglossaryentry{uh}{name=\ensuremath{\boldsymbol{u}_{t-h:t}},description={input history at time \gls{t}},user1={\ensuremath{\boldsymbol{u}}}}
\newglossaryentry{dh}{name=\ensuremath{\boldsymbol{d}_{t-h:t}},description={disturbance history at time \gls{t}},user1={\ensuremath{\boldsymbol{d}}}}
\newglossaryentry{y}{name=\ensuremath{y_t},description={output at time \gls{t}},user1={\ensuremath{y}}}
\newglossaryentry{z}{name={\ensuremath{z_t}},description={nominal state at \gls{t}},user1={\ensuremath{z}}}
\newglossaryentry{e}{name={\ensuremath{e_t}},description={stochastic error term},user1={\ensuremath{e}}}
\newglossaryentry{ey}{name={\ensuremath{e^y_t}},description={error term on the output}}
\newglossaryentry{eu}{name={\ensuremath{e^u_t}},description={error term on the input}}
\newglossaryentry{v}{name={\ensuremath{v_t}},description={nominal input at \gls{t}},user1={\ensuremath{v}}}

\newglossaryentry{zset}{name={\ensuremath{\mathcal{Z}_t}},description={feasible set for \gls{z}}}
\newglossaryentry{eset}{name={\ensuremath{\mathcal{E}_t}},description={error set for \gls{e}}}
\newglossaryentry{esety}{name={\ensuremath{\mathcal{E}^y_t}},description={error set for $\gls{e}^y$}}
\newglossaryentry{esetu}{name={\ensuremath{\mathcal{E}^u_t}},description={error set for $\gls{e}^u$}}

\newglossaryentry{l}{name=\ensuremath{l(\gls{x}, \gls{u})},description={cost at time \gls{t}}}
\newglossaryentry{L}{name=\ensuremath{L(\gls{x}[T])},description={terminal cost at time $T$}}
\newglossaryentry{f}{name=\ensuremath{\boldsymbol{f}(\gls{x}, \gls{u}, \gls{d})},description={dynamics modeling the state transition from \gls{t} to \gls{t}$+1$}}

\newglossaryentry{h}{name={\ensuremath{h}}, description={measurement history}}
\newglossaryentry{m}{name={\ensuremath{m}}, description={available data points}}

\newglossaryentry{db}{name={\ensuremath{\mathcal{F}}},description={dataset for \gls{adp} algorithms}}
\newglossaryentry{Y}{name={\ensuremath{Y}},description={Hankel matrix of measurements \gls{y}}}
\newglossaryentry{U}{name={\ensuremath{\Upsilon}},description={Hankel matrix of applied control inputs \gls{u}}}
\newglossaryentry{D}{name={\ensuremath{D}},description={Hankel matrix of disturbances \gls{d}}}

\newglossaryentry{nx}{name={\ensuremath{n_x}},description={state dimension}}
\newglossaryentry{nu}{name={\ensuremath{n_u}},description={input dimension}}
\newglossaryentry{nd}{name={\ensuremath{n_d}},description={disturbance dimension}}
\newglossaryentry{ny}{name={\ensuremath{n_y}}, description={output dimension}}

\newglossaryentry{G}{name={\ensuremath{G}},description={state space matrices}}
\newglossaryentry{omega}{name={\ensuremath{\Omega}},description={concatenation of \gls{Y}, \gls{U} and \gls{D}}}

\newglossaryentry{A}{name={\ensuremath{A}},description={system matrix}}
\newglossaryentry{B}{name={\ensuremath{B}},description={input matrix}}
\newglossaryentry{E}{name={\ensuremath{E}},description={disturbance matrix}}
\newglossaryentry{C}{name={\ensuremath{C}},description={output matrix}}

\newglossaryentry{Atilde}{name={\ensuremath{\tilde{A}}},description={system matrix in the original space}}
\newglossaryentry{Btilde}{name={\ensuremath{\tilde{B}}},description={input matrix in the original space}}
\newglossaryentry{Etilde}{name={\ensuremath{\tilde{E}}},description={disturbance matrix in the original space}}

\newglossaryentry{Xspace}{name={\ensuremath{\mathcal{X}^{\gls{nx}}}},description={state space},user1={\ensuremath{\mathcal{X}}}}
\newglossaryentry{Uspace}{name={\ensuremath{\mathcal{U}^{\gls{nu}}}},description={control input space},user1={\ensuremath{\mathcal{U}}}}
\newglossaryentry{Dspace}{name={\ensuremath{\mathcal{D}^{\gls{nd}}}},description={disturbance space},user1={\ensuremath{\mathcal{D}}}}

\newglossaryentry{mu}{name={\ensuremath{\mu}},description={penalty term for model regularization}}

\newglossaryentry{p}{name={\ensuremath{p_{t}}},description={energy prices at \gls{t}},user1={\ensuremath{p_t^T}}}
\newglossaryentry{pel}{name={\ensuremath{p_{t,el}}},description={electricity prices at \gls{t}}}
\newglossaryentry{pth}{name={\ensuremath{p_{t,th}}},description={heat prices at \gls{t}}}

\newglossaryentry{Q}{name={\ensuremath{Q_t}},description={state-input value function},user1={\ensuremath{Q}}}
\newglossaryentry{J}{name={\ensuremath{J_t}},description={state value function},user1={\ensuremath{J}}}

\newglossaryentry{mdot}{name={\ensuremath{\dot{m}}},description={mass flow rate},user1={\ensuremath{\dot{m}}}}
\newglossaryentry{dpr}{name={\ensuremath{dp}},description={differential pressure}}

\newglossaryentry{sy}{name={\ensuremath{s^{y}_t}},description={slack variable on output constraints},user1={\ensuremath{s^{y}}}}
\newglossaryentry{su}{name={\ensuremath{s^{u}_t}},description={slack variable on input constraints},user1={\ensuremath{s^{u}}}}

\newglossaryentry{w}{name={\ensuremath{w}},description={optimization penalty term},user1={\ensuremath{w^T}}}

\newglossaryentry{prob_x}{name={\ensuremath{p_x}},description={probability on state constraint violation}}
\newglossaryentry{prob_u}{name={\ensuremath{p_u}},description={probability of input constraint violation}}

\newglossaryentry{K}{name={\ensuremath{K_t}},description={time-varying feedback controller}}

\newglossaryentry{dmean}{name={\ensuremath{\overline{d}_t}},description={mean disturbance},user1={\ensuremath{\overline{d}}}}
\newglossaryentry{dzero}{name={\ensuremath{\tilde{d}_t}},user1={\ensuremath{\tilde{d}}},description={zero mean disturbance}}
\newglossaryentry{dcov}{name={\ensuremath{\Sigma^{\glsuseri{dzero}}_t}},description={covariance matrix for \gls{dzero}}}

\newglossaryentry{xcov}{name={\ensuremath{\Sigma^{\gls{x}}_t}},description={covariance matrix for \gls{x}},user1={\ensuremath{\Sigma^{\gls{x}}}}}

\newglossaryentry{normaldensity}{name={\ensuremath{\phi}}, description={density of normal distribution}}

\newglossaryentry{probtilde_x}{name={\ensuremath{\tilde{p}_x}}, description={violation probability for single constraints}}

\newglossaryentry{freq}{name={\ensuremath{f}},description={compressor frequency}}

\newglossaryentry{tsupply}{name={\ensuremath{T_s}},description={supply temperature}}

\newglossaryentry{tset}{name={\ensuremath{T_{set}}},description={setpoint for \gls{tsupply}}}

\newglossaryentry{ethout}{name={\ensuremath{E_{th, out}}},description={heat pump's thermal energy output [kWh]}}
\newglossaryentry{ethin}{name={\ensuremath{E_{th, in}}},description={heat pump's thermal energy input [kWh]}}
\newglossaryentry{eel}{name={\ensuremath{E_{el}}},description={heat pump's electrical energy input [kWh]}}

\newglossaryentry{trvi}{name={\ensuremath{TRV_i}},description={thermostatic radiator valve for room $i$},text={\ensuremath{TRV_{\glsargi}}}}

\newglossaryentry{tmin}{name={\ensuremath{\underline{y}_t}},description={minimum room temperatures at \gls{t}},user1={\ensuremath{\underline{y}}}}
\newglossaryentry{tmax}{name={\ensuremath{\overline{y}_t}},description={maximum room temperatures at \gls{t}},user1={\ensuremath{\overline{y}}}}

\newglossaryentry{umin}{name={\ensuremath{\underline{u}_t}},description={minimum input at \gls{t}},user1={\ensuremath{\underline{u}}}}
\newglossaryentry{umax}{name={\ensuremath{\overline{u}_t}},description={maximum input at \gls{t}},user1={\ensuremath{\overline{u}}}}

\newglossaryentry{ti}{name={\ensuremath{T_i}},description={room temperature},user1={\ensuremath{T}}}

\newglossaryentry{ghi}{name={\ensuremath{\psi}},description={global horizontal irradiation}}
\newglossaryentry{tdb}{name={\ensuremath{T_{db}}},description={dry bulb temperature}}

\newglossaryentry{cop}{
name={\ensuremath{COP_{hp}}},
description={heat pump's coefficient of performance}
}

\newglossaryentry{numberofrooms}{
name={\ensuremath{N_r}},
description={number of rooms}
}

\newglossaryentry{numberofpreds}{
name={\ensuremath{N_p}},
description={number of predictions}
}

\newglossaryentry{predictedtemp}{
name={\ensuremath{\hat{T}_i}},
description={predicted \gls{ti}},
user1={\ensuremath{\hat{T}}}
}

\newglossaryentry{predictedTsupply}{
name={\ensuremath{\hat{T}_s}},
description={predicted \gls{tsupply}}
}

\newglossaryentry{eocc}{
name={\ensuremath{E_{occ}}},
description={internal gains from occupant behavior}
}

\newglossaryentry{templb}{
name={\ensuremath{\underline{T}_i}},
description={lower temperature constraint for room $i$}
}

\newglossaryentry{tempub}{
name={\ensuremath{\overline{T}_i}},
description={upper temperature constraint}
}

\newglossaryentry{eps}{name={\ensuremath{\epsilon}},description={exploration parameter}}

\newglossaryentry{cp}{name={\ensuremath{c_p}},description={specific heat capacity}}

\newglossaryentry{Uf}{name={\ensuremath{g_j(\gls{x})}},description={observables of state \gls{x}}, user1={\ensuremath{g_j}}}
\newglossaryentry{Ufext}{name={\ensuremath{g_j(\gls{x}, \gls{u}, \gls{d})}},description={observables of state \gls{x}, inputs \gls{u} and disturbances \gls{d}}, user1={\ensuremath{g_j}}}
\newglossaryentry{F}{name={\ensuremath{\boldsymbol{F}^{\Delta t}(\gls{x})}},description={flow map of a dynamical system with timestep $\Delta t$},user1={\ensuremath{\boldsymbol{F}^{\Delta t}}}}
\newglossaryentry{Koop}{name={\ensuremath{\mathcal{K}}},description={Koopman operator of a dynamical system},user1={\ensuremath{\mathcal{K}}}}
\newglossaryentry{kapprox}{name={\ensuremath{K}},description={approximation of the Koopman operator of a dynamical system},user1={\ensuremath{K}}}
\newglossaryentry{khat}{name={\ensuremath{\hat{K}}},description={approximation of the Koopman operator of a dynamical system in a subspace},user1={\ensuremath{\hat{K}}}}
\newglossaryentry{eigf}{name={\ensuremath{\varphi_i(\gls{x})}},description={eigenfunctions of the Koopman operator for autonomous systems},user1={\ensuremath{\varphi_i}}}
\newglossaryentry{eigfext}{name={\ensuremath{\varphi_i(\gls{x}, \gls{u}, \gls{d})}},description={eigenfunctions of the Koopman operator for systems with inputs and disturbances},user1={\ensuremath{\varphi}}}
\newglossaryentry{eigfbasis}{name={\ensuremath{\theta_k(\gls{x}, \gls{u}, \gls{d})}},description={dictionary function of eigenfunction space},user1={\ensuremath{\theta}}}
\newglossaryentry{eigv}{name={\ensuremath{\lambda_i}},description={eigenvalues of the Koopman operator}}
\newglossaryentry{kmode}{name={\ensuremath{\boldsymbol{v}_i}},description={Koopman modes of the Koopman operator}}

\newglossaryentry{embed}{name={\ensuremath{\boldsymbol{\Gamma}(\gls{x})}},description={embedding of the state \gls{x}},user1={\ensuremath{\boldsymbol{\Gamma}}}}
\newacronym{ehp}{EHP}{Electric Heat Pump}
\newacronym{dr}{DR}{Demand Response}
\newacronym{res}{RES}{Renewable Energy Sources}
\newacronym{bess}{BESS}{Battery Energy System}
\newacronym{soc}{SoC}{State of Charge}
\newacronym{hvac}{HVAC}{Heating, Ventilation and Air Conditioning}
\newacronym{dam}{DAM}{Day-Ahead Market}
\newacronym[longplural={Balance Responsible Parties}]{brp}{BRP}{Balance Responsible Party}
\newacronym{ev}{EV}{Electric Vehicle}
\newacronym{pv}{PV}{Photovoltaics}
\newacronym{ewh}{EWH}{Electric Water Heater}
\newacronym{fcr}{FCR}{Frequency containment reserve}
\newacronym{dhw}{DHW}{Domestic Hot Water}
\newacronym{vpp}{VPP}{Virtual Power Plant}
\newacronym{mpc}{MPC}{Model Predictive Control}
\newacronym{dp}{DP}{Dynamic Programming}
\newacronym{adp}{ADP}{Approximate Dynamic Programming}
\newacronym{rhc}{RHC}{Receding Horizon Control}
\newacronym{lti}{LTI}{Linear Time Invariant}
\newacronym{pi}{PI-controller}{Proportional Integral controller}
\newacronym{brl}{BRL}{Batch Reinforcement Learning}
\newacronym{smpc}{SMPC}{Stochastic Model Predictive Control}
\newacronym{scmpc}{ScMPC}{Scenario-based Model Predictive Control}
\newacronym{rmpc}{RMPC}{Robust Model Predictive Control}
\newacronym{socp}{SOCP}{Second Order Cone Problem}
\newacronym{fqi}{FQI}{Fitted q-iteration}
\newacronym{ifqi}{IFQI}{Informed FQI}
\newacronym[plural=TCLs]{tcl}{TCL}{Thermostatically Controlled Load}
\newacronym{cnn}{CNN}{Convolutional Neural Network}
\newacronym{mlp}{MLP}{Multilayer Perceptron}
\newacronym{lstm}{LSTM}{Long Short-Term Memory}
\newacronym{mdp}{MDP}{Markov Decision Process}
\newacronym{pomdp}{POMDP}{Partially Observable Markov Decision Process}
\newacronym{tou}{ToU}{Time of Use}
\newacronym{rbc}{RBC}{Rule Based Controller}
\newacronym{hc}{HC}{Hysteresis Controller}
\newacronym{rl}{RL}{Reinforcement Learning}
\newacronym{dhn}{DHN}{District Heating Network}
\newacronym{ltdhn}{LTDHN}{Low Temperature District Heating Network}
\newacronym{ideas}{IDEAS}{Integrated District Energy Assessment Simulations}
\newacronym{dmd}{DMD}{Dynamic Mode Decomposition}
\newacronym{dmdc}{DMDc}{Dynamic Mode Decomposition with Control}
\newacronym{edmd}{EDMD}{Extended Dynamic Mode Decomposition}
\newacronym{svd}{SVD}{Singular Value Decomposition}
\newacronym{pod}{POD}{Proper Orthogonal Decomposition}
\newacronym{trv}{TRV}{Thermostatic Radiator Valve}
\newacronym{gfs}{GFS}{Global Forecast System}
\newacronym{mad}{MAE}{Mean Absolute Error}
\newacronym{rmse}{RMSE}{Root Mean Square Error}
\newacronym{piv}{PI}{Prediction Interval}
\newacronym{ml}{ML}{Maximum Likelihood}
\newacronym{admm}{ADMM}{Alternating Direction Method of Multipliers}

\title{A concise introduction to Koopman operator theory and the Extended Dynamic Mode Decomposition}
\author[1,2]{Christophe Patyn \thanks{Corresponding author email address: \texttt{christophepatyn@protonmail.com}}}
\author[1,2]{Geert Deconinck}
\affil[1]{Department of Electrical Engineering, KU Leuven, Leuven, Belgium}
\affil[2]{EnergyVille, Genk, Belgium}

\begin{document}

\maketitle

\begin{abstract}
The framework of Koopman operator theory is discussed along with its connections to \gls{dmd} and (Kernel) \gls{edmd}. This paper provides a succinct overview with consistent notation. The authors hope to provide an exposition that more naturally emphasizes the connections between theory and algorithms which may result in a sense of clarity on the subject.
\end{abstract}

\tableofcontents

\section*{Notation}
Vectors are bold lowercase, as in \gls{x}. Lowercase letters (not bold), such as $c_i$, are scalars. Capital letters are matrices, e.g. $X$, while the dual is indicated by $X^*$. A subscript can indicate the element in a vector, such as $g_1(\gls{x})$, a vector that is a row/column of a matrix, as in $\boldsymbol{v}_i$, or may refer to a time index, as in $\gls{x}$. Superscripts indicate powers of a scalar, matrix or operator. There are two exceptions to this: $X^{-1}$ is the matrix inverse and $X^{\dagger}$ is the pseudoinverse. $c_i \in \mathbb{R}$ indicates $c_i$ is in the set of reals. $\mathcal{M}$ indicates the manifold of the system's state space. $\mathcal{K}$ is an operator. $\langle \cdot \rangle$ is the inner product. $\hat{X}$ is an approximation of matrix $X$.

\section{Koopman operator theory}


Koopman showed that every nonlinear dynamical system has an equivalent infinite-dimensional, but globally linear representation \cite{koopman1931hamiltonian, koopman1932dynamical}. Instead of a nonlinear action on a state, the Koopman operator acts linearly on observables \gls{Uf}. These observables may be linear or nonlinear. They may be the identity acting on the full state or a part of it. As an example, a scalar-valued temperature measurement is a single observable, and so is a vector-valued temperature measurement of e.g. the temperatures of multiple zones in a house. They can be combinations of state elements as well, such as the total energy in the system. The resulting linear representation can be analyzed with all the usual tools applicable to linear systems, hence the allure of Koopman operator theory. 
The following is a summary of Koopman operator theory based on multiple expositions in other works \cite{brunton2021modern,rowley2009spectral,lan2013linearization,proctor2018generalizing}.

\subsection{Koopman decomposition}

Starting from a continuous-time, finite-dimensional, autonomous dynamical system as in \cref{eq:autonomousdynsys}, it is assumed that the state \gls{x} evolves on an \gls{nx}-dimensional smooth manifold $\mathcal{M}$. This essentially means the state space is continuous and differentiable, or in other words well-behaved. For more information on topology and manifolds, see Schuller's lectures on differential geometry \cite{schullerdiffgeom}.

\begin{equation} \label{eq:autonomousdynsys}
	\dot{\glsuseri{x}}_t = \boldsymbol{f}(\gls{x})
\end{equation}

The same well-behavedness is assumed for the flow map $\gls{F}: \mathcal{M} \rightarrow \mathcal{M}$, which describes the discrete-time behavior of the system. It maps solutions at time $t$ to solutions at $t+1$. 

\begin{equation} \label{eq:flow}
	\glsuseri{x}[_{\gls{t}+1}] = \gls{F}
\end{equation}

Scalar-valued observable functions $\gls{Uf}: \mathcal{M} \rightarrow \mathbb{R}$ then map from the state space to the reals. The Koopman operator \gls{Koop} acts on these functions and is defined in \cref{eq:koopman}. It is a globally linear operator that maps a, possibly nonlinear, observable one timestep into the future.

\begin{equation} \label{eq:koopman}
	\gls{Koop} \gls{Uf} = \glsuseri{Uf}[(\gls{F})]
\end{equation}

Since the Koopman operator is linear, it has an eigendecomposition with an infinite number of eigenfunctions:

\begin{equation} \label{eq:koopmaneigen}
	\gls{Koop}\gls{eigf} = \gls{eigv} \gls{eigf}; \; i=1,...,\infty
\end{equation}

These eigenfunctions can be nonlinear. They are the generalization of eigenvectors to nonlinear dynamical systems. They essentially measure important quantities such as the system's energy or entropy content and hence encode structures in the state space such as invariant sets or periodic orbits \cite{budivsic2012geometry,mauroy2013isostables}. Furthermore, comparing \cref{eq:koopman} with \cref{eq:koopmaneigen}, the eigenfunctions are the ideal observables, since they diagonalize the Koopman operator and thus provide a basis for the observables' function space. Therefore, we can write every possible observable \gls{Uf} as a weighted sum of eigenfunctions with the coefficients $c_i$ as the weights:

\begin{equation} \label{eq:observablesexpansion}
	\gls{Uf} = \sum^{\infty}_{i=1} \gls{eigf} c_i
\end{equation}

Stacking \gls{ny} observables into a vector-valued observable $\boldsymbol{g}(\gls{x})$ results in \cref{eq:vectorobservable}. This equation shows how the expansion can be interpreted in two dual ways. On the one hand, $\boldsymbol{g}(\gls{x})$ lies in the span of the nonlinear eigenfunctions of the dynamical system. In the dual interpretation, $\boldsymbol{g}(\gls{x})$ lies in the span of the linear Koopman modes \gls{kmode}. These are linear vectors that depend on the chosen observables \gls{Uf} \cite{surana2016koopman}.

\begin{equation} \label{eq:vectorobservable}
	\boldsymbol{g}(\gls{x}) = \begin{bmatrix}
				g_1(\gls{x}) \\
				g_2(\gls{x}) \\
				g_3(\gls{x}) \\
				\vdots \\
				g_{\gls{ny}}(\gls{x})
				\end{bmatrix} = \sum^{\infty}_{i=1} \gls{eigf} \begin{bmatrix}
				c_{1,i} \\
				c_{2,i} \\
				c_{3,i} \\
				\vdots \\
				c_{\gls{ny},i}
				\end{bmatrix} = \sum^{\infty}_{i=1} \gls{eigf} \gls{kmode}
\end{equation}

Filling this out in the left side of \cref{eq:koopman} and using the vector-valued observable $\boldsymbol{g}(\gls{x})$, the Koopman decomposition is obtained (\cref{eq:koopdecomp}) as a combination of three basic structures, being Koopman eigenvalues, Koopman eigenfunctions and Koopman modes. The eigenfunctions describe time-invariant directions in the space of observables, while the modes describe time-invariant directions in the state space defined by the eigenfunctions. The eigenvalues, just like in linear systems, describe growth rates and oscillations along either of these directions \cite{brunton2021modern}.

\begin{equation} \label{eq:koopdecomp}
	\gls{Koop} \boldsymbol{g}(\gls{x}) = \gls{Koop} \sum^{\infty}_{i=1} \gls{eigf} \gls{kmode} = \sum^{\infty}_{i=1} \gls{Koop} \gls{eigf} \gls{kmode} = \sum^{\infty}_{i=1} \gls{eigv} \gls{eigf} \gls{kmode}
\end{equation}

\subsection{Koopman operator theory for linear autonomous systems}

\begin{equation} \label{eq:linearauton}
	\glsuseri{x}[_{t+1}] = A \gls{x}
\end{equation}

Following the example in Rowley et al. \cite{rowley2009spectral}, consider a linear discrete-time system as in \cref{eq:linearauton}. The eigenvectors of $A$ are $\boldsymbol{v}_i$ with the usual eigenvalue expression $A \boldsymbol{v}_i = \lambda_i \boldsymbol{v}_i$. The left eigenvectors of $A$ are the eigenvectors of the adjoint matrix $A^*$, with $A^* \boldsymbol{w}_i = \lambda_i \boldsymbol{w}_i$. They can be normalized so that the inner product of right and left eigenvectors equals the Kronecker delta, i.e. $\varphi_j(\boldsymbol{v}_k) = \langle \boldsymbol{w}_j, \boldsymbol{v}_k \rangle = \delta_{jk}$, where $\langle \cdot \rangle$ is the inner product. The left eigenvectors can be proven to be eigenfunctions \gls{eigf} of the Koopman operator, with $\gls{eigf} = \langle \boldsymbol{w}_i, \gls{x} \rangle$. And the eigenvalues of $A$ are also the eigenvalues of the Koopman operator \gls{Koop}.

\begin{equation} \label{eq:linexample}
	\gls{Koop} \gls{eigf} = \glsuseri{eigf}[(A \gls{x})] = \langle \boldsymbol{w}_i, A \gls{x} \rangle = \langle A^* \boldsymbol{w}_i, \gls{x} \rangle = \lambda_i \langle \boldsymbol{w}_i, \gls{x} \rangle = \lambda_i \gls{eigf}
\end{equation}

And as in \cref{eq:koopdecomp}, the Koopman decomposition can now be obtained with the eigenvectors of $A$ being the Koopman modes, as long as $A$ has a full set of eigenvectors. Note that even though $A$ may have a finite number of eigenvectors, \gls{Koop} has a countably infinite number of eigenvalues and eigenfunctions since $\lambda_i^k \lambda_j^l$ is also an eigenvalue with eigenfunction $\varphi_i^k(\gls{x})\varphi_j^l(\gls{x})$ \cite{brunton2021modern}.

\begin{equation}
	\sum^{\infty}_{i=1} \gls{Koop} \gls{eigf} \gls{kmode} = \sum^{\infty}_{i=1} \gls{eigv} \gls{eigf} \gls{kmode} = \sum^{\infty}_{i=1} \gls{eigv} \langle \boldsymbol{w}_i, \gls{x} \rangle \gls{kmode}
\end{equation}

The eigenfunctions of the Koopman operator are in this case linear functions since the dynamical system is linear. For nonlinear systems, the eigenfunctions will be nonlinear as well, and are tangent to the eigenvectors of the linearized system at the fixed point. So there is a clear connection between eigenfunctions and eigenvectors. The eigenvalues, on the other hand, will equal the eigenvalues of the linearization of the nonlinear system at the fixed point \cite{lan2013linearization}.

\section{Approximating the Koopman operator in practice}
The Koopman operator is infinite-dimensional, so in order to do computation based on this theory a finite-dimensional approximation will have to be made. Techniques to obtain such approximations mainly find their origin in the reduced-order modeling space. Without going too much into detail, they are strongly rooted in approaches such as the \gls{pod} (more commonly known as principal component analysis) and empirical Galerkin projection \cite{holmes2012turbulence}. The \gls{dmd} algorithm is probably the most popular algorithm to approximate the Koopman operator. Its connection to the Koopman operator is shown here.

\subsection{Dynamic mode decomposition}
The aim of \gls{dmd} is to approximate the Koopman modes \gls{kmode} of the Koopman operator. The vector-valued observables $\boldsymbol{g}(\gls{x})$ at each timestep are organised in a snapshot matrix, as in \cref{eq:snapshota}. Considering that the Koopman operator propagates the observables through time, \cref{eq:snapshotb} is obtained.

\begin{subequations}
	\begin{align}
	X & = \big[ \boldsymbol{g}(\gls{x}), \;  \boldsymbol{g}(\glsuseri{x}[_{t+1}]),  \; \boldsymbol{g}(\glsuseri{x}[_{t+2}]),  \; ... \; \boldsymbol{g}(\glsuseri{x}[_{t+n}]) \big] \label{eq:snapshota} \\ 
	  & = \big[ \boldsymbol{g}(\gls{x}),  \; \gls{Koop}\boldsymbol{g}(\gls{x}),  \; \glsuseri{Koop}[^2]\boldsymbol{g}(\gls{x}),  \; ... \; \glsuseri{Koop}[^n]\boldsymbol{g}(\gls{x}) \big] \label{eq:snapshotb}
	\end{align}
\end{subequations}

\gls{dmd} assumes:

\begin{itemize}
	\item there is a linear mapping, i.e. the Koopman operator, between subsequent snapshots;
	\item as the number of snapshots increases, the column vectors inside the snapshot matrix become linearly dependent.
\end{itemize}

As a result of this, the observable at timestep $t+n$ can be written as a linear combination of the other columns. The following equation arises:

\begin{equation} \label{eq:lincombocomp}
	X' = X C
\end{equation}

The companion matrix $C$ is defined as:

\begin{equation}
	C = \begin{bmatrix}
		0 &         &        &   & c_1 \\
		1 & 0 	   &        &    & c_2 \\
		  & \ddots  & \ddots &   & \vdots \\
		  &         & 1      & 0 & c_{n} \\
		  &         &        & 1 & c_{n+1}
	\end{bmatrix}
\end{equation}

The eigenvectors of $C$ are given by the rows of Vandermonde matrix $T$, which contain a sequence of powers of the eigenvalues of $C$, as in \cref{eq:vandermonde}.

\begin{equation} \label{eq:vandermonde}
	T = \begin{bmatrix}
		1       & \lambda_1     & \lambda_1^2     & \ldots & \lambda_1^n \\
		1       & \lambda_2     & \lambda_2^2     & \ldots & \lambda_2^n \\
		\vdots  & \vdots        & \vdots         & \ddots & \vdots \\
		1       & \lambda_n & \lambda_n^2 & \ldots & \lambda_n^n \\
		1       & \lambda_{n+1}     & \lambda_{n+1}^2     & \ldots & \lambda_{n+1}^n

	\end{bmatrix}
\end{equation}

The consequence is that \cref{eq:lincombocomp} can be further decomposed:
\begin{equation}
	X' = X C = X T^{-1} \Lambda T = V \Lambda T
\end{equation}

Since $X T^{-1} = V$, the decomposition of the original $X$ matrix is defined as $X = V T$, which results in a decomposition of the observables (the columns of $X$) as:

\begin{equation}
	\boldsymbol{g}(\glsuseri{x}[_{t+m}]) = \sum_{i=1}^{n+1} \lambda_i^{m} \gls{kmode}
\end{equation}

Comparing this to \cref{eq:koopdecomp}, it is clear that the eigenvalues of $C$ are approximations of the Koopman eigenvalues, while the modes \gls{kmode} are approximations of the Koopman modes scaled by the constant values of the eigenfunctions \gls{eigf} \cite{rowley2009spectral,schmid2010dynamic}.

In order to compute this approximation, an Arnoldi method may be used, as is described in Schmid et al. \cite{schmid2010dynamic}. However, this algorithm is ill-conditioned. Instead, the \gls{svd} of the snapshot matrix is computed \cite{tu2013dynamic}, i.e. $X = U \Sigma W^T$, to obtain the following equation:

\begin{subequations}  \label{eq:svddmd}
	\begin{align}
		& X' = \gls{Koop} X = \gls{Koop} U \Sigma W^T \\
		\Longleftrightarrow \; & U^T X' W \Sigma^{-1} = U^T \gls{Koop} U = \gls{khat}
	\end{align}
\end{subequations}

\gls{khat} is the approximation of the Koopman operator. The Koopman modes can then be approximated as $\gls{kmode} = \frac{1}{\gls{eigv}} X' W \Sigma^{-1} \boldsymbol{w}_i$, with $\boldsymbol{w}_i$ the i'th eigenvector of $\hat{\gls{Koop}}$ (for a proof, see \cite{tu2013dynamic}). In reduced-order modeling, the $U$ matrix is said to be the \gls{pod} basis, which holds the main directions of correlation between the scalar-valued observables \gls{Uf}. This can be interpreted as a projection of the Koopman operator onto the \gls{pod} basis, represented by the column space of $U$ which links back to the assumption of linear dependency of the columns in the snapshot matrix. The \gls{svd} enables accounting for rank-deficiency by removing the columns in $U$ which correspond to close-to-zero singular values, resulting in a more robust approximation algorithm \cite{schmid2010dynamic,tu2013dynamic}. 

Locally, \gls{dmd} allows to construct a linear model near steady states where the dynamics are only weakly nonlinear. Global validity is theoretically achievable if this linear model maps observables that lie in the span of the eigenfunctions of the system \cite{cenedese2022data}.

\subsection{Koopman operator for non-autonomous systems}
So far, the Koopman operator has been described in the context of autonomous systems in the form of \cref{eq:autonomousdynsys} (for continuous time), and \cref{eq:flow} (for discrete time). These do not consider control inputs and disturbances, which, in practice, do need to be taken into account. The eigenfunctions \gls{eigf} in \cref{eq:koopmaneigen} are only defined as a function of the state \gls{x}. However, they can be expanded to include functions of control inputs \gls{u} and disturbances \gls{d} and any cross terms, resulting in the eigenfunctions \gls{eigfext}. Proctor et al. \cite{proctor2018generalizing} describe a number of approaches to define these. The approach relevant to this work is the one where the inputs and disturbances are constant during a timestep and are generated from an exogenous forcing term (such as a controller or the ambient temperature). As a result, the Koopman operator is defined as in \cref{eq:koopmanext}.

\begin{equation} \label{eq:koopmanext}
	\gls{Koop} \boldsymbol{g}(\gls{x}, \gls{u}, \gls{d}) = \boldsymbol{g}(\glsuseri{F}(\gls{x}, \gls{u}, \gls{d}), \gls{u}, \gls{d})
\end{equation}

With this new observable space, the Koopman decompostion now comes in terms of these new eigenfunctions, as in \cref{eq:koopdecompext}.

\begin{equation} \label{eq:koopdecompext}
	\gls{Koop} \boldsymbol{g}(\gls{x}, \gls{u}, \gls{d}) = \gls{Koop} \sum^{\infty}_{i=1} \gls{eigfext} \gls{kmode} = \sum^{\infty}_{i=1} \gls{Koop} \gls{eigfext} \gls{kmode} = \sum^{\infty}_{i=1} \gls{eigv} \gls{eigfext} \gls{kmode}
\end{equation}

\subsection{Koopman operator for partially observable systems}
Now that inputs and disturbances are accounted for, partial observability of the state is addressed. Usually, the state is not observed in full. Rather, measurements of a few quantities of the state space (i.e. observable functions on the state space) are available. Fitting \gls{khat} on this small number of observables, will usually not result in any significantly accurate model. Taken's embedding theorem \cite{rand_detecting_1981} (and the later work by Sauer et al. \cite{sauer1991embedology}) provides a solution. This theorem states that for a smooth dynamical system $\gls{F}: \mathbb{R}^{\gls{nx}} \rightarrow \mathbb{R}^{\gls{nx}}$ and a smooth observable $\gls{Uf}: \mathbb{R}^{\gls{nx}} \rightarrow \mathbb{R}^{\gls{ny}}$, the map $\gls{embed}: \mathbb{R}^{\gls{nx}} \rightarrow \mathbb{R}^{2\gls{nx}+1}$ is an embedding, with \gls{embed} defined in \cref{eq:embed}. 

\begin{equation} \label{eq:embed}
	\gls{embed} = (\boldsymbol{g}(\gls{x}), \; \boldsymbol{g}(\boldsymbol{F}^{\Delta t}(\gls{x})), \; ..., \; \boldsymbol{g}(\boldsymbol{F}^{2\gls{nx}\Delta t}(\gls{x})))
\end{equation}

Since an embedding is a smooth, bijective map from the state space to the space of delay-embedded measurements, its inverse also exists \cite{huke2006embedding}. Therefore, a new dynamical system is obtained in \cref{eq:embeddynsys} as a composition of the embedding and the original dynamical system, with the history of observables defined as $\glsuseri{yh}[_{t:t+h}] = (\boldsymbol{g}(\gls{x}), \; \boldsymbol{g}(\boldsymbol{F}^{\Delta t}(\gls{x})), \; ..., \; \boldsymbol{g}(\boldsymbol{F}^{2\gls{nx}\Delta t}(\gls{x})))$, with $h=2\gls{nx}+1$.

\begin{equation} \label{eq:embeddynsys}
	\boldsymbol{F}^{\Delta t}_T(\glsuseri{yh}[_{t:t+h}]) = \glsuseri{embed}( \glsuseri{F}(\boldsymbol{\Gamma}^{-1}(\glsuseri{yh}[_{t:t+h}])))
\end{equation}

The implication of this theorem is that taking a history of observables, and defining this as the state of the delay-embedded dynamical system $\boldsymbol{F}^{\Delta}_T(\cdot)$, results in this new system having similar characteristics to the original system \cite{huke2006embedding}. Taken's theorem was later generalized to input-output systems \cite{casdagli1992dynamical}, resulting in a delay-embedded system with control inputs and disturbances as in equation \cref{eq:delayinputsdists}, with:

\begin{equation}
 \glsuseri{yh}[_{t:t+h}] = (\boldsymbol{g}(\gls{x}, \gls{u}, \gls{d}), \; \boldsymbol{g}(\boldsymbol{F}^{\Delta t}(\gls{x}, \glsuseri{u}[_{t:t+1}], \glsuseri{d}[_{t:t+1}])), \; ..., \; \boldsymbol{g}(\boldsymbol{F}^{2\gls{nx}\Delta t}(\gls{x}, \glsuseri{u}[_{t:t+h}], \glsuseri{d}[_{t:t+h}])))
\end{equation}

\begin{equation} \label{eq:delayinputsdists}
	\glsuseri{yh}[_{t+1:t+h+1}] = \boldsymbol{F}^{\Delta t}_T(\glsuseri{yh}[_{t:t+h}], \gls{u}, \gls{d})
\end{equation}

$\glsuseri{u}[_{t:t+i}]$ and $\glsuseri{d}[_{t:t+i}]$ are histories of inputs and disturbances (not observables) from timestep $t$ until $t+i$. \gls{dmd} can be applied to the delay-embedded dynamical system to approximate its Koopman operator. The snapshot matrix then consists of columns of a history of observables $\glsuseri{yh}[_{t:t+h}]$. Multiple authors have proposed such an approach \cite{arbabi2017ergodic,arbabi2018data,brunton2017chaos,hirsh2021structured} for both autonomous and non-autonomous systems with partial observability. One of the more important results, albeit theoretical, is that for an autonomous ergodic system, the delay embedding, in the limit of infinite time, results in the retrieval of the Koopman eigenfunctions and eigenvalues of the original dynamical system \cite{arbabi2017ergodic}. 

\subsection{Extended dynamic mode decomposition} \label{sec:edmd}
In the previous sections, the assumption has been that the observables \gls{Ufext} lie in the subspace spanned by the eigenfunctions \gls{eigfext} \cite{williams2015data}:

\begin{equation} \label{eq:obsinv}
	\boldsymbol{g}(\gls{x}, \gls{u}, \gls{d}) = \sum^n_{i=1} \gls{eigfext} \gls{kmode} 
\end{equation}

Since \gls{eigfext} are eigenfunctions of the Koopman operator by definition, the subspace spanned by the eigenfunctions is invariant under the action of the Koopman operator. Therefore, given the assumption defined by \cref{eq:obsinv}, the space of observables is also invariant under the action of the Koopman operator \cite{kutz2016dynamic,tu2013dynamic,brunton2021modern}. Additionally, the previous section showed that the delay-embedded system $\boldsymbol{F}^{\Delta t}_T(\cdot)$ is a composition of the embedding and the flow map of the original system. Most likely the original system or the embedding, or both, are nonlinear. The delay-embedded system therefore will be nonlinear as well. As a consequence, the linear Koopman operator with its linear modes can only be defined if the eigenfunctions are nonlinear. Therefore, using the \gls{dmd} or \gls{dmdc} algorithm results in an approximation of the eigenfunctions to first order \cite{williams2015data}, which again alludes to its local (rather than global) validity \cite{cenedese2022data}. \gls{edmd} solves this problem. From now on the delay-embedded $\glsuseri{yh}[_{t:t+h}]$ is referred to as the state \gls{x} of the system for notational convenience.

In \gls{edmd} \cite{williams2015data}, a set of real-valued dictionary functions $\{\glsuseri{eigfbasis}[_1], \; \glsuseri{eigfbasis}[_2], \; ..., \glsuseri{eigfbasis}[_k] \}$ is defined. Since $k$ is finite, these are assumed to be a basis for the eigenfunctions of the system, as formalized in \cref{eq:eigfbasis1} and \cref{eq:eigfbasis}. This set has to be rich enough to allow for accurate prediction, but may be restricted to a subspace based on domain knowledge.

\begin{equation} \label{eq:eigfbasis1}
	\gls{eigfext} = \sum_{k=1}^K b_{ik} \gls{eigfbasis}
\end{equation}

or as matrix-vector multiplication, where $B$ contains the coefficients $b_{ik}$ in the dictionary functions basis:

\begin{subequations} \label{eq:eigfbasis}
	\begin{align}
		\boldsymbol{\phi}(\gls{x}, \gls{u}, \gls{d}) & = B \boldsymbol{\theta}(\gls{x}, \gls{u}, \gls{d}) \\
		\text{with } & \boldsymbol{\phi}(\gls{x}, \gls{u}, \gls{d}) = \big[ \glsuseri{eigfext}[_1](\gls{x}, \gls{u}, \gls{d}), \; 
										 \glsuseri{eigfext}[_2](\gls{x}, \gls{u}, \gls{d}), \; 
										 ..., \;
										 \glsuseri{eigfext}[_n](\gls{x}, \gls{u}, \gls{d})  \big]^T \\
		\text{ }  & \boldsymbol{\theta}(\gls{x}, \gls{u}, \gls{d}) = \big[ \glsuseri{eigfbasis}[_1](\gls{x}, \gls{u}, \gls{d}), \; 
										 \glsuseri{eigfbasis}[_2](\gls{x}, \gls{u}, \gls{d}), \; 
										 ..., \;
										 \glsuseri{eigfbasis}[_k](\gls{x}, \gls{u}, \gls{d})  \big]^T
	\end{align}
\end{subequations}


The dictionary functions are used to construct the new snapshot matrix:

\begin{equation} \label{eq:thetacoll}
	X = \big[ \boldsymbol{\theta}(\gls{x}, \gls{u}, \gls{d}), \; 
			    \boldsymbol{\theta}(\glsuseri{x}[_{t+1}], \glsuseri{u}[_{t+1}], \glsuseri{d}[_{t+1}]), \;
			     ..., \; 
			     \boldsymbol{\theta}(\glsuseri{x}[_{t+n}], \glsuseri{u}[_{t+n}], \glsuseri{d}[_{t+n}])  \big]
\end{equation}

Following from the equation $X' = \gls{khat} X$, the Koopman operator can be obtained by linear regression, i.e. $\gls{khat} = X'X^{\dagger}$. The same \gls{svd} approach as in  \cref{eq:svddmd} can be used to reduce computational cost and add regularization. \gls{khat} maps the dictionary functions at $t$ to those at $t+1$:  $\boldsymbol{\theta}(\glsuseri{x}[_{t+1}], \glsuseri{u}[_{t}], \glsuseri{d}[_{t}]) = \gls{khat} \boldsymbol{\theta}(\gls{x}, \gls{u}, \gls{d})$. Diagonalizing \gls{khat} results in \cref{eq:diagkhat_theta}, where $P$ are the eigenvectors of \gls{khat}, and $\Lambda$ the diagonal eigenvalue matrix.

\begin{equation} \label{eq:diagkhat_theta}
	\gls{khat} \boldsymbol{\theta}(\gls{x}, \gls{u}, \gls{d}) = P \Lambda P^{-1} \boldsymbol{\theta}(\gls{x}, \gls{u}, \gls{d})
\end{equation}

Combining \cref{eq:eigfbasis,eq:diagkhat_theta}, it becomes clear that the rows of matrix $B$ are the eigenvectors of the dual of the Koopman operator approximation since $P$ diagonalizes \gls{khat}. Therefore, to obtain the eigenfunctions, first \gls{khat} is identified with the \gls{edmd} algorithm. Then its dual eigenvectors, i.e. the rows of $P^{-1}$, are used in \cref{eq:eigfbasis} as the $B$ matrix.

To approximate the Koopman modes, the observables are assumed to lie in the span of the dictionary functions, i.e. $\boldsymbol{g}(\gls{x}, \gls{u}, \gls{d}) = D \boldsymbol{\theta}(\gls{x}, \gls{u}, \gls{d})$ where $D$ holds the coefficients of this expansion in the basis defined by the dictionary functions. If this is not the case, then $D$ can be approximated by linear regression. Again, the accuracy of the approximation depends on the chosen function basis. Combining this equation with \cref{eq:eigfbasis} results in \cref{eq:deriv_modes}, where $V$ is the approximation of the Koopman modes that results from the \gls{edmd} algorithm. 

\begin{equation} \label{eq:deriv_modes}
	\boldsymbol{g}(\gls{x}, \gls{u}, \gls{d}) = D \boldsymbol{\theta}(\gls{x}, \gls{u}, \gls{d}) =  D B^{-1} \boldsymbol{\phi}(\gls{x}, \gls{u}, \gls{d}) = V \boldsymbol{\phi}(\gls{x}, \gls{u}, \gls{d})
\end{equation}

\subsection{Kernel trick for EDMD}

As \cref{eq:thetacoll} shows, the dictionary functions are collected in a matrix $X \in \mathbb{R}^{k \times n}$ and the approximation can be computed as $\gls{khat} = X'X^{\dagger}$. The issue with this approach is that the set of $k$ dictionary functions can be so large that $\gls{khat} \in \mathbb{R}^{k \times k}$ will be computationally intractable to compute. For this reason, the kernel trick was proposed as a computationally tractable way to implement \gls{edmd} \cite{williams2014kernel}. In the original approximation \gls{khat}, the pseudo-inverse of $X$ can be computed by means of the \gls{svd} of $X$, as in $X^{\dagger} = W \Sigma^{-1} U^T$. Therefore, \gls{khat} lies in the rowspace of $U^T$ and every eigenvector of \gls{khat} can be written as a combination of the columns of $U \in \mathbb{R}^{k \times n}$. \Cref{eq:koopapproxedmd} uses this subspace to approximate the Koopman operator as $\glsuseri{khat}[_U] \in \mathbb{R}^{n \times n}$.

\begin{subequations} \label{eq:koopapproxedmd}
	\begin{align}
		\gls{eigv} \boldsymbol{\hat{v}}_i     & = \gls{khat} \boldsymbol{\hat{v}}_i \\
		\gls{eigv} U \boldsymbol{v}_i & = X'X^{\dagger} U \boldsymbol{v}_i \\
							  & = X^{-T} X^T X' X^{\dagger} U \boldsymbol{v}_i \\
							  & = U \Sigma^{-1} W^T (X^T X') W \Sigma^{-1} U^T U \boldsymbol{v}_i \\
							  & = U \big[(\Sigma^{-1} W^T) A ( W \Sigma^{-1} ) \big] \boldsymbol{v}_i \label{eq:kapproxedmd} \\ 
							  & = U \glsuseri{khat}[_U] \boldsymbol{v}_i \\
		\gls{eigv} \boldsymbol{v}_i & = \glsuseri{khat}[_U] \boldsymbol{v}_i
	\end{align}
\end{subequations}

To compute $A$, an inner product of the columns of $X$ and $X'$ is taken. Since these columns contain the dictionary functions defined earlier, every element of $A$ can be written as an inner product:
 
\begin{equation}
	A_{ij}=\boldsymbol{\theta}^T(\glsuseri{yh}[_{i-h:i}], \glsuseri{uh}[_{i-h:i}], \glsuseri{dh}[_{i-h:i}]) \boldsymbol{\theta}(\glsuseri{yh}[_{j-h+1:j+1}], \glsuseri{uh}[_{j-h:j}], \glsuseri{dh}[_{j-h:j}])
\end{equation}

The kernel trick defines a kernel function $k(\boldsymbol{x}, \boldsymbol{y})$ that allows to compute this inner product in the feature space, rather than in the $k$-dimensional state space of \gls{khat}. Multiple different kernels can be used which represent a different set of dictionary functions. For example, a polynomial kernel $k(\boldsymbol{x}, \boldsymbol{y})=(1 + \boldsymbol{x}^T\boldsymbol{y})^{\alpha}$ defines a set of dictionary functions that can represent all polynomials up to and including degree $\alpha$. A radial basis kernel $k(\boldsymbol{x}, \boldsymbol{y}) = e^{-\frac{|| \boldsymbol{x} - \boldsymbol{y} ||}{\sigma^2}}$ represents a polynomial basis of infinite degree. Without the kernel trick, computing the matrix $A$ would require to explicitly compute the dictionary functions, which for example with a polynomial dictionary would already be intractable for low powers, especially considering the history of terms that have to be taken into account in \gls{yh}. What remains is to obtain $W$ and $\Sigma^{-1}$. Using the \gls{svd} of $X$, \cref{eq:xtx} shows how $W$ and $\Sigma^{-1}$ can be derived from the eigenvalue decomposition of a matrix $G$. The same kernel trick can be used to compute every $G_{ij}$ \cite{williams2014kernel}.

\begin{equation}  \label{eq:xtx}
	G = X^T X = W \Sigma U^T U \Sigma W^T = W \Sigma^2 W^T
\end{equation}

To obtain the diagonalized system, the eigenfunctions can be used to project the data. The eigenfunctions are represented by \cref{eq:eigfuncs}, where $k(\gls{x}, \boldsymbol{x}_i)$ is the kernel function of the current state \gls{x} and column $i$ in the data matrix $X$.

\begin{equation} \label{eq:eigfuncs}
	\gls{eigf} = \boldsymbol{\hat{v}}_i^T \Sigma^{-1} Q^T \begin{bmatrix} k(\gls{x}, \boldsymbol{x}_1) \\ \vdots \\ k(\gls{x}, \boldsymbol{x}_n)  \end{bmatrix}
\end{equation}

The vector-valued observable $\boldsymbol{g}(\gls{x})$ can be obtained with \cref{eq:obsvcomp}. Since it is assumed that $\boldsymbol{g}(\gls{x})$ lies in the span of the eigenfunctions, it is a linear function. In this equation, the matrix $V$ contains the eigenvectors of $\glsuseri{khat}[_U]$ in its columns. For more information on the observable and eigenfunctions, refer again to Williams et al. \cite{williams2014kernel}.

\begin{equation} \label{eq:obsvcomp}
	\boldsymbol{g}(\gls{x}) = X Q \Sigma^{-1} V^{-T}
\end{equation}

\section{Conclusion}
In an attempt to provide a concise overview, Koopman operator theory and the theoretical connections with the main algorithms to compute an approximation of the Koopman operator was discussed in this paper. The work is based on information scattered throughout multiple papers on the subject and was presented with notational consistency, hoping to provide clarity on the subject.

\bibliographystyle{ieeetr}
\bibliography{refs}

\end{document}